\documentclass[reqno]{amsart}
\usepackage{hyperref}

\usepackage[centertags]{amsmath}
\usepackage{amsfonts}
\usepackage{amssymb}
\usepackage{amsthm}
\usepackage{newlfont}
\usepackage{fancyhdr}

\makeatletter
\@namedef{subjclassname@2010}{%
  \textup{2020} Mathematics Subject Classification}
\makeatother

\fancypagestyle{plain}{\fancyfoot[C]{\vspace*{2\baselineskip}\thepage}}%Change number to move page number
\theoremstyle{plain}
\newtheorem{theo}{Theorem}[section]
\newtheorem{thmx}{Theorem}[section]
\newtheorem{lem}{Lemma}[section]
\newtheorem{corol}[theo]{corollary}

\newtheorem*{conj*}{Denjoy's Conjecture}

\begin{document}
\title[%EJDE-2012/177
\hfil Higher Order Complex differential equations]
{Fractional order derivative characterizations of Besov-Morrey type spaces with applications}

\author[ \hfil %EJDE-2012/177\hfilneg
]
{Lu Chen, Mingjin Li, Jianren Long}

%\address{Xiubi Wu \newline
%School of Science, Guizhou University, Guiyang, 550025, P.R. China.}
%\email{basicmath@163.com }

\address{Lu Chen \newline
School of Mathematical Sciences, Guizhou Normal
University, Guiyang, 550025, P.R. China.}
\email{2139044243@qq.com}

\address{Mingjin Li \newline
School of Mathematical Sciences, Guizhou Normal
University, Guiyang, 550025, P.R. China.} \email{997208526@qq.com}

\address{Jianren Long \newline
School of Mathematics and Computer Science, Guizhou  Normal
University, Guiyang,  550001, P.R. China.
Corresponding author.}
\email{longjianren2004@163.com}

%\address{Ke-e Qiu \newline
%School of Mathematics and Computer Science, Guizhou Normal
%Colleage, Guiyang, 550018, P.R. China.}
%\email{qke456@sina.com }

\thanks{This research work is supported by National Natural Science Foundation of China (Grant No. 12261023, 11861023).}
%\thanks{Submitted July 6, 2012. Published October 12, 2012.}

\begin{abstract}
  On the one hand, the fractional order derivative characterization of the Besov-Morrey type space $B_{p}^{K}(s)$ is established by $K$-Carleson measures, and it was also shown that \( f \in B_{p}^{K}(s_1) \Leftrightarrow f^{\left(\frac{s_2 - s_1}{p}\right)} \in B_{p}^{K}(s_2) \), which extended the results of Sun et al. on the fractional derivative of Morrey type space. On the other hand, some sufficient conditions for the growth of solutions to linear complex differential equations have been obtained by using $n$th derivative criterion.
\end{abstract}

\maketitle
\numberwithin{equation}{section}
\newtheorem{theorem}{Theorem}[section]
\newtheorem{lemma}[theorem]{Lemma}
\newtheorem{definition}[theorem]{Definition}
\newtheorem{example}[theorem]{Example}
\newtheorem{remark}[theorem]{Remark}
\allowdisplaybreaks

\section{Introduction}

Let $\mathbb{D}$ be the open unit disc in the complex plane and $H(\mathbb{D})$ be the set of all analytic function in $\mathbb{D}$.  For $0< p<\infty,$ the Hardy space $H^{p}$ consists of all functions $f\in H(\mathbb{D})$ such that $$\|f\|_{H^p}^p:=\sup_{0<r<1}\frac1{2\pi}\int_0^{2\pi}|f(re^{i\theta})|^p d\theta<\infty.$$
The classical Hardy space $H^p$ is $H_{0}^p$. For $0<\alpha<\infty $, the $H_\alpha^{\infty}$ spaces consists of all functions $f\in H(\mathbb{D})$ such that $$\|f\|_{H_\alpha^\infty}:=\sup\limits_{z\in\mathbb{D}}|f(z)|(1-|z|^2)^\alpha<\infty.$$
 
Throughout this paper, let $K:[0, \infty) \to [0, \infty)$ be a nondecreasing and right-continous function, not identically zero.
We assume that $K$ satisfies:
 \begin{equation}\label{eq:ad}
  \int_0^1 \frac{\varphi_K(x)}{x} \, dx < \infty
 \end{equation}
and
\begin{equation}\label{eq:ny}
\int_1^\infty\frac{\varphi_K(x)}{x^{1+\sigma}}dx<\infty,\quad\sigma>0,
  \end{equation}where
 $$\varphi_K(x)=\sup\limits_{0<t\leq1}\frac{K(tx)}{K(t)},\quad0<x<\infty.$$

 It is obvious that the function \( K(t) = t^q \) satisfies \eqref{eq:ad} for all \( q \in (0, \infty) \) and satisfies \eqref{eq:ny} whenever \( q\in (0,\sigma) \). There are more weight function results in \cite{i} and \cite{j}.

 For \( 0 \leq p < \infty \), the space \( Q_p \) consists of analytic functions \( f \) on \( \mathbb{D} \) for which

\[
\|f\|_{Q_p}^2 = \sup_{a \in \mathbb{D}} \int |f'(z)|^2 \bigg(g(a,z)\bigg)^p dA(z) < \infty,
\]
where $dA$ is an area measure on $\mathbb{D}$ normalized so that $A(\mathbb{D})$ = 1, $g(a,z) = \log \left| \frac{1-\bar{a}z}{a-z} \right|$.

The \( Q_K \) space consists of analytic functions \( f \) in \( \mathbb{D} \) satisfying$$
    \|f\|_K = \left( \sup_{a \in \mathbb{D}} \int_\mathbb{D} |f'(z)|^2 K(g(z,a)) \, dA(z) \right)^{1/2} < \infty,$$
in the case $K(t) = t^p$, $0 < p < \infty$, the space $Q_K$ gives $Q_p$ space.

 Wu and Xie \cite{a} characterized functions in Morrey spaces using $q$-Carleson measures. Furthermore, they established a simple relationship between $Q$ spaces and Morrey spaces, that is, $Q$ spaces can be viewed as fractional integrals of Morrey spaces, which implies that Morrey spaces correspond to fractional derivatives of $Q$ spaces. The authors proved that the equivalence $f \in \mathcal{Q}_p \Leftrightarrow f^{\left(\frac{1-p}{2}\right)} \in \mathcal{L}^{2,p}$ holds(where \(0 < p \leq 1\) and \(q = \frac{1-p}{2}\)). Moreover, \( g \in \mathcal{L}^{2,p} \) if and only if the measure \( |\nabla g(z)|^2 y \, dx \, dy \) is a bounded \( p \)-Carleson measure. Wulan and Zhu \cite{ww} characterizes the $Q_{K}$ space using higher-order derivatives. Wulan and Zhou \cite{b} characterize $Q_{K}$ space by K-Carleson measure and showed a connection between $H_{K}^{2}$ and $Q_{K}$ by using the fractional order derivative. that is, let \( K \) satisfy the conditions \eqref{eq:ad} and \eqref{eq:ny}, if \( \alpha > 1/2 \), then \( f \in Q_K \) if and only if \( \left| f^{(\alpha)}(z) \right|^2 (1 - |z|^2)^{2(\alpha - 1)} dA(z) \) is a \( K \)-Carleson measure. At this time, their results also indicate that if \( f \in Q_K \), exist  \( 0 < q < p \)(where  $0 < p < \infty$), such that \(f^{\left(\frac{1-q}{2}\right)} \in H^2_K \). The study of higher-order (including fractional order) derivative characterizations in various function spaces has drawn considerable interest among researchers. See, \cite{ku,d,qt,c} for related works.

 Recently, Sun and Wulan \cite{c} obtaied the fractional order derivative characterization of Dirichlet-Morrey type spaces, and the specific results are as follows.
 \begin{thmx}\label{thm:gg}
Suppose \( K \) satisfies (1.2) for some \( 0 < \sigma < 2 \) and \(\frac{\sigma}{2} \leq s < \infty\). Then \( f \in D_K^s \) if and only if  
\begin{align*}
  |f^{(t)}(z)|^2(1 - |z|^2)^{2t-2+s}dA(z), \quad
\max\left\{0, \frac{1-s}{2}, \frac{\sigma-s}{2}\right\} < t < \infty
\end{align*}
is a \( K \)-Carleson measure.
\end{thmx}
 Theorem \ref{thm:gg} actually provides the relationship between $D_{K}^{s_1}$ and $D_{K}^{s_2}$, that is
  \( f \in D_{K}^{s_1} \Leftrightarrow f^{\left(\frac{s_2 - s_1}{2}\right)} \in D_{K}^{s_2} \).

 Let $0<s<1<p<\infty$. A function $f\in H(\mathbb{D})$ belongs to the Besov type space $B_{p}(s)$ if 
$$\|f\|_{B_{p}(s)}^p=|f(0)|^p+\int_\mathbb{D}|f^{\prime}(z)|^p(1-|z|^2)^{p-2+s} dA(z)<\infty,$$
where $dA(z)$ denotes the normalized area measure on $\mathbb{D}$.

In \cite{e}, Qu and Zhu introduced a new Besov-Morrey-type space $B_{p}^{K}(s)$ in their study of Volterra integral operators, defined as follows.

Let $0<s<1<p<\infty$. The Besov-Morrey type space $B_{p}^{K}(s)$ is the space of all $f\in H(\mathbb{D})$ such that 
 $$\|f\|_{B_{p}^{K}(s)}^p=|f(0)|^p+\sup_{a\in\mathbb{D}}\frac{(1-|a|^2)^s}{K(1-|a|^2)}\|f\circ\varphi_a-f(a)\|_{B_{p}(s)}^p<\infty.$$
 where $\varphi_{a}(z)=\frac{a-z}{1-\bar{a}z}$ is a M\"obius transformation of  $\mathbb{D}$. If $ K(t) = t^{s \lambda}$ and $ 0 < s, \lambda < 1 $, then $B_p^K(s) = B_p^{\lambda}(s)$. and if p=2, then  $B_p^K(s) = \mathcal{D}_K^s$. 

Based on Theorem \ref{thm:gg}, Sun and Wulan \cite{c} characterized equivalent norms for the $B_p^K(s)$ space in the case p=2. A natural question arises whether similar results hold for $p\neq 2$, which is precisely the focus of the present paper.

Before presenting the main results, we also need to review some notations.

For an arc \(I \subset \partial \mathbb{D}\) is an interval, \(|I| = \frac{1}{2\pi} \int_{I} |dz|\) denotes the normalized length of \(I\), the Carleson square is denoted by
\[
S(I) = \{re^{it} : e^{it} \in I, 1 - |I| \leq r < 1\}.
\]
 A positive Borel measure \(\mu\) on \(\mathbb{D}\) is called a \(K\)-Carleson measure provided
\[
\sup_{I \subset \partial \mathbb{D}} \frac{\mu(S(I))}{K(|I|)} < \infty.
\]
If \(K(t) = t^s\) (\(0 < s < \infty\)), then \(\mu\) is a Carleson measure. Specially, \(\mu\) is called the classical Carleson measure whenever \(s = 1\).

The paper is organized as follows. In section 1, relevant background and review of some symbols are provided. In section  2, the fractional order derivative of the Besov-Morrey type space is obtained. In section 3, this fractional order derivative is applied to linear complex differential equations, and the sufficient condition for all analytical solutions of equation (1.1) to belong to the Besov-Morrey type space as coefficients is obtained.

In this paper, $A\lesssim B$ means that there exists a constant $C = C(\cdot)$ (it can be slightly different each time it appears) such that $A\leq CB$. If both $A\lesssim B$ and $B\lesssim A$, then we write $A\approx B$.

\section{Fractional order derivative characterization of Besov-Morrey type space $B_{p}^{K}(s)$}

In this section, we will characterize the fractional order derivatives of Besov-Morrey type space. 
For fixed \( b > 1 \), define the \(\alpha\)-order derivative as follows:

\[f^{(\alpha)}(z) = \frac{\Gamma(b + \alpha)}{\Gamma(b)} \int_{D} \frac{(1 - |w|^2)^{b-1}}{(1 - \overline{w}z)^{b+\alpha}} \overline{w}^{[\alpha-1]} f'(w) \, dA(w), \quad b + \alpha > 0,\]

where \(\Gamma\) is the Gamma function and \([\alpha]\) denotes the smallest integer which is larger than or equal to \(\alpha\). Since

\[(z^n)^{(\alpha)} = 
\begin{cases} 
\frac{\Gamma(b+n+\alpha-1-[\alpha-1]) \Gamma(n+1)}{\Gamma(b+n)\Gamma(n-[\alpha-1])} z^{n-1-[\alpha-1]}, & n \geq [\alpha-1] + 1, \\
0, & n < [\alpha-1] + 1,
\end{cases}\]
If \(\alpha \in N^+\), then \(f^{(\alpha)}\) is just the derivative of order \(n\) of \(f\).

Now, we give the main results of this section.

\begin{theo}\label{thm:a}
Suppose $K$ satisfies the (1.2) for some $\sigma\in (0,2s)$, $0<s<1$, $\max\{1, 1+\sigma-s\}<p<\infty $, and $f\in H(\mathbb{D})$. Then $f\in B_{s}^{K}(s)$ if and only if
\begin{align*}
\left|f^{(t)}(z)\right|^p (1-|z|^2)^{pt-2+s} \, dA(z), \quad
\frac{2-s}{p}+\frac{p}{p-1}< t < \infty
\end{align*}
is a K-carleson measure.
\end{theo}

\begin{remark} For the case of $p=2$ in Theorem 2.1, this is the result in Theorem \ref{thm:gg}.
\end{remark}

Theorem \ref{thm:a} actually gives a relationship between two spaces $B_{p}^{K}(s_1)$ and $B_{p}^{K}(s_2)$ as follows.

\begin{corol} $K$ satisfies the (1.2) for some $\sigma\in (0,2s)$, $0<s<1$, $\max\{1, 1+\sigma-s\}<p<\infty $, For \(0 < s_1, s_2 < 1\) with \( s_1 < s_2 + p \), the following statements are equivalent:
\begin{enumerate}
    \item \( f \in B_{p}^{K}(s_1) \).
    \item \( f^{\left(\frac{s_2 - s_1}{p}\right)} \in B_{p}^{K}(s_2) \).
\end{enumerate}
\end{corol}

Before proving the theorem, we will need some auxiliary lemmas.

\begin{lem}\textnormal{\cite[Theorem2.1]{c}}\label{lem:zkl}
Suppose $K$ satisfies condition (1.2) for some $\sigma \in (0, 2)$, and $\mu$ is a positive Borel measure on $\mathbb{D}$. Then $\mu$ is a $K$-Carleson measure if and only if
\[
\sup_{a \in \mathbb{D}} \frac{1}{K(1 - |a|^2)} \int_{\mathbb{D}} \left( \frac{1 - |a|^2}{|1 - \overline{a}z|} \right)^q d\mu(z) < \infty, \quad \sigma \leq q < \infty.
\]
\end{lem}

\begin{lem}\textnormal{\cite[Remark2.1(i)]{c}}\label{lem:mbq}
Suppose \( 0<t \leq r < \infty \), \( K \) satisfies condition (1.2) for some \(\sigma > 0\). Then
\[
\frac{K(r)}{K(t)} \leq \left( \frac{r}{t} \right)^{\sigma-c} \leq \left( \frac{r}{t} \right)^{\sigma}
\]
holds for all sufficiently small constant \( c < \sigma \).
\end{lem}

\begin{lem}\textnormal{\cite[Proposition 5]{e}}\label{lem:hjj}
Let \( 0 < s < 1 < p < \infty \), \( f \in H(\mathbb{D}) \), and suppose \( K \) satisfies condition (1.2) for some \( \sigma \in (0, 2s) \). Then \( f \in B_p^K(s) \) if and only if
\[
\sup_{I \subset \partial \mathbb{D}} \frac{1}{K(|I|)} \int_{S(I)} |f'(z)|^p (1 - |z|^2)^{p-2+s} \,dA(z) < \infty.
\]
\end{lem}

The following integral estimates play a fundamental role in both the study of function space operators on the unit disk and the proof of Lemma \ref{lem:xz}.

\begin{lem}\textnormal{\cite[Lemma 3.10]{f}}\label{lem:xzg}
 Suppose \( z \in \mathbb{D} \), \( c \) is real, \( t > -1 \). Then
\begin{align*}
I_{c,t}(z) &= \int_{\mathbb{D}} \frac{(1-|w|^2)^t}{|1-z\overline{w}|^{2+t+c}} \, dA(w), \\
 &\asymp 
\begin{cases} 
1, & c < 0, \\ 
\dfrac{1}{(1-|z|^2)^c}, & c > 0,\ |z| \to 1^-, \\ 
\log\dfrac{1}{1-|z|^2}, & c = 0,\ |z| \to 1^-. 
\end{cases}
\end{align*}
\end{lem}

We need the following technical lemmas.

\begin{lem}\label{lem:xz}
   Suppose $K$ satisfies the (1.2) for some $\sigma\in (0,2)$, $0<s<1<p<\infty$,
\begin{align*}
 \frac{2-s}{p}+\frac{p}{p-1}< \alpha < \infty, \quad 2+p+\frac{s-2}{p} < b < \infty.
\end{align*}
For $f \in L^1(\mathbb{D}, dA)$, define 
$$Tf(z) = \int_{\mathbb{D}} \frac{(1 - |w|^2)^{b-1}}{|1 - \overline{w}z|^{\alpha + b}} f(w) \, dA(w), \quad z \in \mathbb{D}.$$
If $|f(z)|^p (1 - |z|^2)^{p-2+s} \, dA(z)$
is a K-Carleson measure, then
$$|Tf(z)|^p (1 - |z|^2)^{p\alpha-2+s} \, dA(z)$$
is also a $K$-Carleson measure.
\end{lem}

\begin{proof}
For given \( p \), $s$, let
\[
d\mu_f(z) = |f(z)|^p(1-|z|^2)^{p-2+s}dA(z),
\]
\[
d\mu_{T_f}(z) = |Tf(z)|^p (1 - |z|^2)^{p\alpha-2+s} dA(z),
\]
\[
\|\mu\|_K = \sup_{I \subset \partial \mathbb{D}} \frac{\mu(S(I))}{K(|I|)}.
\]

Next we will show that \(\|\mu_f\|_K < \infty\) implies \(\|\mu_{Tf}\|_K < \infty\). For \( I \subset \partial \mathbb{D} \), suppose the centre of \( 2^n I \subset \partial\mathbb{D} \) is same with \( I \), both of their length are \( 2^n |I| \), where \( n \in \mathbb{N}^+ \). Then
\begin{align*}
\mu_{Tf}(S(I)) &= \int_{S(I)} |Tf(z)|^p (1 - |z|^2)^{p\alpha - 2 + s} dA(z) 
\\&\leq \int_{S(I)} \left( \int_{S(2I)} + \int_{\mathbb{D}\setminus S(2I)} \frac{|f(w)| (1 - |w|^2)^{b-1}}{|1-\overline{w}z|^{\alpha + b}} dA(w) \right)^p (1-|z|^2)^{p\alpha - 2 + s}dA(z)
\\&\lesssim E_1 + E_2,
\end{align*}
where
\[
E_1 = \int_{S(I)} \left( \int_{S(2I)} \frac{|f(w)| (1 - |w|^2)^{b-1}}{|1 - \overline{w}z|^{\alpha + b}} dA(w) \right)^p (1 - |z|^2)^{p\alpha - 2 + s}dA(z)
\]
\[
E_2 = \int_{S(I)} \left( \int_{\mathbb{D}\setminus S(2I)} \frac{|f(w)| \left( 1 - |w|^2 \right)^{b-1}}{|1 - \overline{w}z|^{\alpha + b}} \, dA(w) \right)^p \left( 1 - |z|^2 \right)^{p\alpha - 2 + s} \, dA(z).
\]
For \( E_1 \), we consider
\[
H(z,w) = \frac{(1 - |z|^2)^{\alpha-\frac{2}{p} + \frac{s}{p}}(1 - |w|^2)^{b-2+\frac{2}{p}-\frac{s}{p}}}{|1 - \overline{w}z|^{\alpha + b}}
\]
and the integral operator

\[
T_H g(z) = \int_{\mathbb{D}} |g(w)| H(z,w) \, dA(w), \quad g \in L^p(\mathbb{D}).
\]
where \[
g(w) = \left(1 - |w|^2\right)^{1-\frac{2-s}{p}} |f(w)| \chi_{S(2I)}(w)
\]and \(\chi_{S(2I)}\) is the characteristic function of \(S(2I)\).
We obtain that
\[
\int_{\mathbb{D}} H(z, w)(1 - |z|^2)^p \, dA(w) \lesssim (1 - |w|^2)^p
\]
and
\[
\int_{\mathbb{D}} H(z, w)(1 - |w|^2)^\frac{p}{p-1} \, dA(z) \lesssim (1 - |z|^2)^\frac{p}{p-1}.
\]

Accordingly\textnormal{\cite[Theorem 3.6]{f}}, the operator \( T_H \) is bounded from \( L^p(\mathbb{D}) \) to \( L^p(\mathbb{D}) \).
We have
\[
 E_1 \leq \int_{\mathbb{D}} \left( \int_{\mathbb{D}} |g(w)| H(z,w) \, dA(w) \right)^p dA(z) \leq \int_{\mathbb{D}} |g(z)|^p \, dA(z) \leq \| \mu_f \|_K K(|I|). 
\]
To handle \( E_2 \), we note that both
\begin{equation}\label{eq:w}
  2^n|I| \leq |1 - \overline{w}z|, \quad z \in S(I), \quad w \in S(2^{n+1}I) \setminus S(2^nI),
\end{equation}
and
\begin{equation}\label{eq:e}
\int_{S(2^nI)} (1 - |z|^2)^{\tau-2} \, dA(z) \lesssim (2^n|I|)^\tau, \quad \tau > 1,
\end{equation}
hold for all \( n \in \mathbb{N} \). Combining \eqref{eq:w}, \eqref{eq:e}, \ref{lem:mbq} and H\"older inequality, we get
\begin{align*}
E_2 &\lesssim \int_{S(I)} \left( \sum_{n=1}^{\infty} (2^n|I|)^{-(\alpha + b)} \int_{S(2^{n+1}I)} \frac{|f(w)|}{(1 - |w|^2)^{1-b}} \, dA(w) \right)^p (1 - |z|^2)^{p\alpha - 2 + s}dA(z)
\\&\lesssim |I|^{p\alpha + s} \left( \sum_{n=1}^{\infty} (2^n|I|)^{-(\alpha + b)} \int_{S(2^{n+1}I)} \frac{|f(w)|}{(1 - |w|^2)^{1-b}} \, dA(w) \right)^p.
\\&\lesssim |I|^{p\alpha+s} \left[ \sum_{n=1}^{\infty} (2^n |I|)^{(-\alpha+b)} \left( \int_{S(2^{n+1} I)} |f(w)|^p (1 - |w|^2)^{p-2+s}\, dA(w) \right)^{\frac{1}{p}} 
\right.
\\&\text{ }\text{ }\left. \cdot \left( \int_{S(2^{n+1} I)} (1 - |w|^2)^{\frac{p(b-2)+2-s}{p-1}} \, dA(w) \right)^{\frac{p-1}{p}} \right]^p
\\&\lesssim |I|^{p\alpha+s} \left[ \sum_{n=1}^{\infty} (2^n |I|)^{-(\alpha+b)} \cdot \mu_{f}(S(2^{n+1} I))^{\frac{1}{p}} \cdot (2^{n+1} |I|)^{b-\frac{s}{p}} \right]^p.
 \\&\lesssim |I|^{p\alpha+s} \left[ \sum_{n=1}^{\infty} (2^n |I|)^{-\alpha-\frac{s}{p}} \cdot K^{\frac{1}{p}}(2^{n+1} |I|) \right]^p \cdot \|\mu_{f}\|_K
\\&\lesssim \left[ \sum_{n=1}^{\infty} 2^{n(-\alpha-\frac{s}{p})} \cdot 2^{\frac{(n+1)\sigma}{p}} \right]^p \cdot K(|I|) \cdot \|\mu_{f}\|_K.
\\&\lesssim \left[ \sum_{n=1}^{\infty} 2^{-n(\alpha+\frac{s-\sigma}{p})} \right]^p \cdot K(|I|) \cdot \|\mu_{f}\|_K.
\\&\lesssim  K(|I|) \cdot \|\mu_{f}\|_K.
\end{align*}
The foregoing estimates on \( E_1 \) and \( E_2 \) give \( \mu_{Tf} \) is a \( K \)-Carleson measure.
\end{proof}

Now we give the proof of Theorem \ref{thm:a}.

\noindent\textbf{Proof of Theorem 2.1.}
 Let \( f \in B_{p}^{K}(s) \). By Lemma \ref{lem:hjj}, we know that \( |f'(z)|^{p}(1-|z|^{2})^{p-2+s}dA(z) \) is a \( K \)-Carleson measure. Since

\[
|f^{(t)}(z)| \leq \frac{\Gamma(b+t)}{\pi\Gamma(b)} \int_{\mathbb{D}} \frac{(1-|w|^{2})^{b-1}|f'(w)|}{|1-\overline{w}z|^{b+t}} dA(w), \quad 2+p+\frac{s-2}{p} < b < \infty,
\]

we obtain by Lemma \ref{lem:xz} that

\[
|f^{(t)}(z)|^{p}(1-|z|^{2})^{pt-2+s}dA(z)
\]

is a \( K \)-Carleson measure.

Conversely, suppose \( |f^{(t)}(z)|^{p}(1-|z|^{2})^{pt-2+s}dA(z) \) is a \( K \)-Carleson measure. Now we consider the function:

\[
g(z) = \frac{\Gamma(b+1)z^{[t-1]}}{\pi\Gamma(b+t-1)} \int_{\mathbb{D}} \frac{(1-|w|^{2})^{b+t-2}}{(1-\overline{w}z)^{b+1}} f^{(t)}(w) dA(w).
\]

It follows from Lemma \ref{lem:xz} that \( |g(z)|^{p}(1-|z|^{2})^{p-2+s}dA(z) \) is a \( K \)-Carleson measure. Next, we will consider the relationship between $g$ and \( f' \).

Let
\[ 
f(z) = \sum_{j=0}^{\infty} a_j z^j
\]
and
\[ 
f^{(t)}(z) = \sum_{j=0}^{\infty} a_{j,t} z^j, \quad z \in \mathbb{D},
\]
where
\[
a_{j,t} = a_{j+m+1} \frac{\Gamma(j+b+t)\Gamma(j+m+2)}{\Gamma(j+1)\Gamma(j+m+1+b)}, \quad j \in \mathbb{N}, \, m = [t-1].
\]

Obviously, \( m \geq 0 \). Rewrite

\[
g(z) = \sum_{j=0}^{\infty} \frac{B(j+b+1,m)}{B(j+1,m)} (j+m+1)a_{j+m+1}z^{j+m},
\]
where \( B(\cdot, \cdot) \) is the Beta function. Note that

\[
f'(z) = g(z) + \sum_{j=0}^{\infty} (j + m + 1) a_{j+m+1} z^{j+m} + \sum_{j=1}^{m} j a_j z^{j-1} - g(z).
\]

Let
\[
s_m(z) = \sum_{j=1}^{m} j a_j z^{j-1}
\]

and

\[
h(z) = \sum_{j=0}^{\infty} \left( 1 - \frac{B(j + b + 1, m)}{B(j + 1, m)} \right) a_{j+m+1} z^{j+m+1},
\]

we obtain

\begin{equation}
f'(z) = g(z) + s_m(z) + h'(z).
\end{equation}

\textbf{Case 1:} When \( m = 0 \), we can easily obtain \( f' = g \). Therefore, \( |f'(z)|^p (1 - |z|^2)^{p-2+s} dA(z) \) is a \( K \)-Carleson measure.

\textbf{Case 2:} When \( m > 0 \), since \( |g(z)|^p (1 - |z|^2)^{p-2+s} dA(z) \) is a \( K \)-Carleson measure and \( s_m \) is a polynomial. We need only to show
that \( |h'(z)|^p (1 - |z|^2)^{p-2+s} dA(z) \) is also a K-Carleson measure. By Lemma \ref{lem:hjj}, we get
\begin{align*}
\frac{1}{K(|I|)} \int_{S(I)} |h'(z)|^p (1 - |z|^2)^{p-2+s} dA(z) &\leq \frac{1}{|I|^\sigma} \int_{S(I)} |h'(z)|^p (1 - |z|^2)^{p-2+s} dA(z)
\\ &\leq \int_{\mathbb{D}} |h'(z)|^p (1 - |z|^2)^{p-2+s-\sigma} dA(z) 
\\&\leq \|h\|_{p,p-2+s-\sigma,0}^{p}.
\end{align*}

By\cite[Theorem 5.5]{g}, we have

\begin{equation}\label{eq:wyb}
\|f\|_{p,q,0}^p = \sum_{k=0}^{\infty} (n_k)^{p-q-1} |a_k|^p, \quad 0< p < \infty, \quad -1 < q < \infty.
\end{equation}

 and a standard estimate

\[
0 < 1 - \frac{B(j + b + 1, m)}{B(j + 1, m)} \leq \frac{(b + 1)m}{j + m + 1}, \quad j \in \mathbb{N}.
\]
Accordingly,

\[
\|h\|_{_{p,p-2+s-\sigma,0}}^{p} \approx \sum_{j=0}^{\infty} \left( 1 - \frac{B(j + b + 1, m)}{B(j + 1, m)} \right)^{p} (j + m + 1)^{1-s+\sigma} |a_{j+m+1}|^{p}
\]

\[
\lesssim \sum_{j=0}^{\infty} \frac{|a_{j+m+1}|^{p}}{(j + m + 1)^{p+s-1-\sigma}}.
\]

In addition, by Lemma \ref{lem:zkl} for any \( q \), \(\sigma \leq q < \infty\),

\begin{align*}
\infty &> \sup_{I \subseteq \partial \mathbb{D}} \frac{1}{K(|I|)} \int_{S(I)} |f^{(t)}(z)|^{p} (1 - |z|^{2})^{pt-2+s} dA(z)
\\&\approx \sup_{a \in \mathbb{D}} \frac{1}{K(1 - |a|^{2})} \int_{\mathbb{D}} \left( \frac{1 - |a|^{2}}{|1 - \overline{a} z|} \right)^{q} |f^{(t)}(z)|^{p} (1 - |z|^{2})^{pt-2+s} dA(z).
\\&\gtrsim \int_{\mathbb{D}} |f^{(t)}(z)|^2 (1 - |z|^2)^{pt-2+s} \, dA(z) 
\\&\approx \sum_{j=0}^{\infty} \frac{|a_{j+m+1}|^2}{(j+m+1)^{-s}}.
\end{align*}

Here, the latter part holds by using \eqref{eq:wyb} and the Stirling's formula:

\[
\frac{\Gamma(n+c)}{n!} \approx n^{c-1}, \quad c > 0.
\]

Therefore,

\[
\|h\|_{p,p-2+s-\sigma,0}^p \lesssim \sum_{j=0}^{\infty} \frac{|a_{j+m+1}|^p}{(j+m+1)^{-s}} < \infty,
\]
which is the desired result.
\qed

The following corollary is obtained through a simple calculation combining Theorem \ref{thm:a} and Lemma \ref{lem:zkl}, so we omit the proof.
\begin{corol}\label{cor:qr}
Suppose $K$ satisfies the (1.2) for some $\sigma\in (0,s)$, $0<s<1$, $1+\sigma-s<p<\infty$. let $n$ be a positive integer and $f\in H(\mathbb{D})$. Then $f\in B_{p}^{K}(s)\Leftrightarrow \sup\limits_{a \in \mathbb{D}}I_n<\infty$,
where \[I_n:=
\sup_{a \in \mathbb{D}} \frac{(1 - |a|^2)}{K(1 - |a|^2)} \int_{\mathbb{D}} |f^{(n)}(z)|^p (1 - |z|^{np-4+s})(1 - |\sigma_a(z)|^2)^2 \, dA(z) < \infty.
\]
\end{corol}
\section{Applications}
The growth of solutions to linear complex differential equations is an interesting topic that has been extensively studied by many scholars. They have obtained both fast growing and slowly growing solutions using different approaches. For fast growing solutions, Nevanlinna theory have been successfully applied, see \cite{hk,ll,nu}. Conversely, the study of slowly growing solutions often requires various distinct methods, such as Herold's comparison, Picard's approximations and Gronwall's Lemma, see \cite{pp,k,bb}.

In 1982, Pommerenke \cite{do} studied the complex second-order differential equation
\begin{equation}\label{eq:xq}
  f^{''}+A_0f=0,
\end{equation}
where $A(z)$ is an analytic function in the unit disk $\mathbb{D}$. Sufficient conditions for the coefficient function $A(z)$ such
that all solutions of equation \eqref{eq:xq} belong to Hardy spaces are first found by Carleson measure. Subsequently, numerous results on sufficient conditions for all solutions of higher-order equations to belong to other function spaces have been obtained.

Heittokangas \cite{k} et al. obtained the expression for the solution of non-homogeneous linear complex differential equations
\begin{equation}\label{eq:hh}
  f^{(k)}+A_{k-1}(z)f^{(k-1)}+\cdots+A_1(z)f^{'}+A_0f=A_k(z),
\end{equation}
 and their research on the spatial properties of the solution of complex differential equation \eqref{eq:hh} attracted widespread attention from scholars. Pel\'{a}ez and R\"{a}tty\"{a} \cite{qw} discovered that when the coefficients of equation \eqref{eq:hh}  satisfy certain conditions, the solutions of this equation belong to the weighted Bergman space \( A^p_\omega \). Hausko et al. \cite{n} obtained coefficient conditions under which the solutions of equation \eqref{eq:hh}  belong to the weighted Hardy space \( H^\infty_\omega \). More research results have been obtained in \cite{o,p,l,m,jj}.

We mainly obtained the following results.
\begin{theo}\label{thm:b}
Suppose $K$ satisfies the (1.2) for some $\sigma\in (0,s)$, and $n$ be a positive integer. Let $A_j\in H(\mathbb{D})$, $j=0,1,\cdots,n$, $0<s<1$, $1+\sigma-s<p<\infty $, and $f\in H(\mathbb{D})$. If
\[
M_0 := \sup_{a \in \mathbb{D}} \frac{(1 - |a|^2)^2}{K(1 - |a|^2)} \int_{\mathbb{D}} |A_0(z)|^p(1 - |z|^2)^{np-4} (1 - |\sigma_a(z)|^2)^2 \, dA(z),
\]  
and  
\[
M_1 := \sum_{j=1}^{n-1} \|A_j\|_{H_{n-j}^\infty},\] 
are both sufficiently small positive constants,
\[
M_2 := \sup_{a \in \mathbb{D}} \frac{(1 - |a|^2)^2}{K(1 - |a|^2)} \int_{\mathbb{D}} |A_n(z)|^p K(1 - |z|^2)(1 - |z|^2)^{np-4+s} (1 - |\sigma_a(z)|^2) \, dA(z) < \infty,
\]  
Then all solutions of \eqref{eq:hh} belong to $B_{p}^{K}(s)$.

\end{theo}

The following theorem is obtained through the integration condition to obtain the sufficient condition for the coefficients of equation \eqref{eq:hh} when the solution belongs to $B_{p}^{K}(s)$.

\begin{theo}\label{thm:fly}
Suppose $K$ satisfies the (1.2) for some $\sigma\in (0,s)$, and $n$ be a positive integer. Let $A_j\in H(\mathbb{D})$, $j=0,1,\cdots,n$, $0<s<1$, $1+\sigma-s<p<\infty $, and $f\in H(\mathbb{D})$. If
\[
N_0 := \sup_{a \in \mathbb{D}} \frac{(1 - |a|^2)^2}{K(1 - |a|^2)} \int_{\mathbb{D}} \left| \int_0^z \int_0^{\xi_1} \cdots \int_0^{\xi_{n-2}} |A_0(\xi_{n-1})|^p \frac{K(1 - |\xi_{n-1}|^2)}{(1 - |\xi_{n-1}|^2)^{s}} \, d\xi_{n-1} \cdots d\xi_1 \right|\]
\[
\cdot (1 - |z|^2)^{p-4+s}(1 - |\sigma_a(z)|^2)^2 \, dA(z),
\]
and
\[
N_1 := \sup_{a \in \mathbb{D}} \frac{(1 - |a|^2)^2}{K(1 - |a|^2)} \int_{\mathbb{D}} \bigg( \sum_{m=1}^{n-1} \int_0^z \int_0^{\xi_1} \cdots \int_0^{\xi_{m-1}} \frac{\left|\sum_{k=1}^m A_{n-k}^{(m-k)}(\xi_m)\right|^p}{(1-|\xi_m|^2)^{p+s}} \,
\]
\[
\cdot K(1 - |\xi_{m}|^2)d\xi_m \cdots d\xi_1 \bigg) (1-|z|^2)^{p-4+s}(1-|\sigma_a(z)|^2)^2 \, dA(z)
\]
are both sufficiently small positive constants,
\[
N_2 := \sup_{a \in \mathbb{D}} \frac{(1 - |a|^2)^2}{K(1 - |a|^2)} \int_{\mathbb{D}} \left| \int_0^z \int_0^{\xi_1} \cdots \int_0^{\xi_{n-2}} A_n(\xi_{n-1}) \, d\xi_{n-1} \cdots d\xi_1 \right|^p \]
\[
\cdot(1-|z|^2)^{p-4+s}(1-|\sigma_a(z)|^2)^2 \, dA(z) < \infty,
\]
Then all solutions of \eqref{eq:hh} belong to $B_{p}^{K}(s)$.
\end{theo}

Before proving the theorem, let's first state some necessary lemmas.
\begin{lem}\textnormal{\cite[P.57]{h}}\label{lem:qx}
Suppose that $a_{n}\geq0$ for $n=1,2,\ldots,N$. Then $$\bigg(\sum_{n=1}^Na_n\bigg)^m\leq\bigg(\sum_{n=1}^Na_n^m\bigg) ,\quad0<m\leq1,$$ and $$\bigg(\sum_{n=1}^Na_n\bigg)^m\leq N^{m-1}\bigg(\sum_{n=1}^Na_n^m\bigg), \quad1\leq m<\infty.$$
\end{lem}
The following lemma can be obtained by combining mathematical induction and
Leibniz’s rule, so the proof is omitted.
\begin{lem}
 Let \( n \) be a natural number. If \( f \in \mathcal{H}(\mathbb{D}) \), then  
\[f^{(n)}(z)h(z) = \sum_{i=0}^{n} (-1)^{i} \binom{n}{i} (f h^{(i)})^{(n-i)}(z), \quad z \in \mathbb{D}.\]
\end{lem}

\begin{lem}\label{lem:cl}
Suppose $K$ satisfies the (1.2) for some $\sigma\in (0,2)$, \(0 < s < 1 < p < \infty\). Let \(n\) be a positive integer. If \(f \in H(\mathbb{D})\), then
\[
|f^{(n)}(a)|^p (1 - |a|^2)^{pn-2+s} \leq \int_{\mathbb{D}} |f^{(n)}(z)|^p (1 - |z|^2)^{pn-4+s} (1 - |\sigma_a(z)|^2)^2 \, dA(z),
\]
where \(a \in \mathbb{D}\).
\end{lem}

\begin{proof}
For any \(a \in \mathbb{D}\), let
\[
E\left(a, \frac{1}{e}\right) = \left\{ z \in \mathbb{D} : |z - a| < \frac{1}{e} (1 - |a|) \right\}.
\]
Obviously, \(E\left(a, \frac{1}{e}\right) \subset D\left(a, \frac{1}{e}\right)\). When \(z \in E\left(a, \frac{1}{e}\right)\), we see
\[
\left( 1 - \frac{1}{e} \right) (1 - |a|) \leq 1 - |z| \leq \left( 1 + \frac{1}{e} \right) (1 - |a|).
\]

By \cite[Proposition 4.13]{f}, we have
\begin{align*}
&\int_{\mathbb{D}} |f^{(n)}(z)|^p (1 - |z|^2)^{pn-4+s} (1 - |\sigma_a(z)|^2)^2 \, dA(z) \\
&\quad \geq \int_{D\left(a, \frac{1}{2}\right)} |f^{(n)}(z)|^p (1 - |z|^2)^{pn-4+s} (1 - |\sigma_a(z)|^2)^2 \, dA(z)  \\
&\quad \geq \int_{E\left(a, \frac{1}{2}\right)} |f^{(n)}(z)|^p (1 - |z|^2)^{pn-4+s} \, dA(z) \\
&\quad \geq (1 - |a|^2)^{pn-4+s} \int_{E\left(a, \frac{1}{2}\right)} |f^{(n)}(z)|^p \, dA(z) \\
&\quad \geq (1 - |a|^2)^{pn-2+s} |f^{(n)}(a)|^p.
\end{align*}

This completes the proof.
\end{proof}

By using Lemma \ref{lem:qx} and Corollary \ref{cor:qr}, we can obtain the following Lemma, which estimates the growth rate of functions in Besov-Morrey type space $B_{p}^{K}(s)$.

\begin{lem}\label{lem:jw}
Suppose $K$ satisfies the (1.2) for some $\sigma\in (0,s)$, $0<s<1$, $1+\sigma-s<p<\infty$. let $n$ be a nonnegative integer and $f\in H(\mathbb{D})$.
\[
|f^{(n)}(z)| \lesssim \left\|f\right\|_{B_{p}^{K}(s)} \bigg(\frac{K(1 - |z|^{2})}{(1 - |z|^{2})^{pn + s}}\bigg)^\frac{1}{p}.
\]
\end{lem}

\noindent\textbf{Proof of Theorem 3.1.}
Suppose that $f$ is an analytic solution of \eqref{eq:hh}, then
\begin{equation}\label{eq:xy}
	f_{r}^{(n)}(z)+\sum_{j=0}^{n-1}B_{j}(z)f_{r}^{(j)}(z)=B_{n}(z), \quad z\in \mathbb{D},
\end{equation}
where $f_{r}(z)=f(rz)$, $B_{j}(z)=B_{j}(z,r)=r^{n-j}A_{j}(r z)$, $0<r<1,j=0,\ldots,n-1$, $B_n(z)=B_n(z,r)=r^n A_j(r z)$. Combining Lemma \ref{lem:qx} and \eqref{eq:xy}, we have
\begin{align*}
\|f_r\|_{B_{p}^{K}(s)}^p &\approx \sup_{a \in \mathbb{D}} \frac{(1 - |a|^2)^2}{K(1 - |a|^2)} \int_{\mathbb{D}} \left| f_r^{(n)}(z) \right|^p (1 - |z|^2)^{np-4+s} (1 - |\sigma_a(z)|^2)^2 \, dA(z)
\\&= \sup_{a \in \mathbb{D}} \frac{(1 - |a|^2)^2}{K(1 - |a|^2)} \int_{\mathbb{D}} \left| B_n(z) - \sum_{j=0}^{n-1} B_j(z) f_r^{(j)}(z) \right|^p
\\&\text{ }\text{ }\text{ }\cdot (1 - |z|^2)^{np-4+s} (1 - |\sigma_a(z)|^2)^2 \, dA(z)
\\&\lesssim \sup_{a \in \mathbb{D}} \frac{(1 - |a|^2)^2}{K(1 - |a|^2)} \int_{\mathbb{D}} \left( |B_n(z)|^p + \sum_{j=0}^{n-1} \left| B_j(z) f_r^{(j)}(z) \right|^p \right)
\\&\text{ }\text{ }\text{ }\cdot (1 - |z|^2)^{np-4+s} (1 - |\sigma_a(z)|^2)^2 \, dA(z)
\\&\lesssim \sup_{a \in \mathbb{D}} \frac{(1 - |a|^2)^2}{K(1 - |a|^2)} \int_{\mathbb{D}} |B_0(z)f_r(z)|^p (1 - |z|^2)^{np-4+s} (1 - |\sigma_a(z)|^2)^2 \, dA(z)
\\&\text{ }\text{ }\text{ }+ \sup_{a \in \mathbb{D}} \frac{(1 - |a|^2)^2}{K(1 - |a|^2)} \left( \sum_{j=1}^{n-1} \int_{\mathbb{D}} |B_j(z)f_r^{(j)}(z)|^p \right)(1 - |z|^2)^{np-4+s} (1 - |\sigma_a(z)|^2)^2 \, dA(z)
\\&\quad+ \sup_{a \in \mathbb{D}} \frac{(1 - |a|^2)^2}{K(1 - |a|^2)} \int_{\mathbb{D}} |B_n(z)|^p (1 - |z|^2)^{np-4+s} (1 - |\sigma_a(z)|^2)^2 \, dA(z)
\\&:= M_{0,1} + M_{1,1} + M_2.
\end{align*}
Using lemma \ref{lem:jw}, we get
\begin{align*}
M_{0,1}&\approx\sup_{a \in \mathbb{D}} \frac{(1 - |a|^2)^2}{K(1 - |a|^2)} \int_{\mathbb{D}} |B_0(z)f_r(z)|^p (1 - |z|^2)^{np-4+s} (1 - |\sigma_a(z)|^2)^2 \, dA(z)
\\&\lesssim\sup_{a \in \mathbb{D}} \frac{(1 - |a|^2)^2}{K(1 - |a|^2)} \int_{\mathbb{D}} |B_0(z)|^p  \left\|f\right\|_{B_{p}^{K}(s)}^p \bigg(\frac{K(1 - |z|^{2})}{(1 - |z|^{2})^{s}}\bigg)(1 - |z|^2)^{np-4+s} (1 - |\sigma_a(z)|^2)^2 \, dA(z)
\\&\lesssim\left\|f\right\|_{B_{p}^{K}(s)}^p\sup_{a \in \mathbb{D}} \frac{(1 - |a|^2)^2}{K(1 - |a|^2)} \int_{\mathbb{D}} |B_0(z)|^p K(1 - |z|^{2})(1 - |z|^2)^{np-4} (1 - |\sigma_a(z)|^2)^2 \, dA(z)
\\&:= \left\|f\right\|_{B_{p}^{K}(s)}^p M_{0}
\end{align*}
We esttimate $M_{1,1}$,
\begin{align*}
M_{1,1}&\lesssim\sup_{a \in \mathbb{D}} \frac{(1 - |a|^2)^2}{K(1 - |a|^2)} \sum_{j=1}^{n-1}\int_{\mathbb{D}} |B_j(z)f_r^{(j)}(z)|^p (1 - |z|^2)^{np-4+s} (1 - |\sigma_a(z)|^2)^2 \, dA(z)
\\&\lesssim\sup_{a \in \mathbb{D}} \frac{(1 - |a|^2)^2}{K(1 - |a|^2)} \sum_{j=1}^{n-1}\int_{\mathbb{D}} |B_j(z)|^p|f_r^{(j)}(z)|^p (1 - |z|^2)^{jp-4+s}(1 - |z|^2)^{np-jp} 
\\&\text{ }\text{ }\text{ }\cdot(1 - |\sigma_a(z)|^2)^2 \, dA(z)
\\&\lesssim\sum_{j=1}^{n-1}\sup_{z \in \mathbb{D}} |B_j(z)|^p (1 - |z|^2)^{np-jp}\sup_{a \in \mathbb{D}} \frac{(1 - |a|^2)^2}{K(1 - |a|^2)} \int_{\mathbb{D}}|f_r^{(j)}(z)|^p (1 - |z|^2)^{jp-4+s} 
\\&\text{ }\text{ }\text{ }\cdot(1 - |\sigma_a(z)|^2)^2 \, dA(z)
\\&:= \left\|f\right\|_{B_{p}^{K}(s)}^p M_{1}^p
\end{align*}
Hence, the conditions of Theorem \ref{thm:b} yield
$$\|f_{r}\|_{B_{p}^{K}}^{p}\lesssim\frac{M_{2}}{1-M_{0}-M_{1}^{p}}<\infty, \quad 0<r<1.$$
This implies $f\in B_{p}^{K}(s)$ by letting $r\to1^{-}$.
\qed

\noindent\textbf{Proof of Theorem 3.2.} Suppose that $f$ is an analytic solution of \eqref{eq:hh}, then
\begin{equation}\label{eq:qc}
	f_r^{(n)}(z)=B_n(z)-\sum_{j=0}^{n-1}B_j(z)f_r^{(j)}(z), \quad z\in \mathbb{D},
\end{equation}
where $f_{r}(z)=f(r z)$, $B_{j}(z)=B_{j}(z,r)=r^{n-j}A_{j}(r z)$, $0<r<1$, $j=0,\ldots,n-1$, $B_n(z)=B_n(z,r)=r^n A_j(rz)$. By applying the identity$$f^{^{\prime}}(z)=\int_0^{z}f^{^{\prime\prime}}(\xi)d\xi+f^{^{\prime}}(0), \quad z\in \mathbb{D},$$
$n-1$ times from 0 to $\xi_{s}$, $s=1,\ldots,n-2$, along any fixed piecewise smooth curve in $\mathbb{D}$. By \eqref{eq:qc}, we have
\begin{align*}
f'_r(z) &= \int_0^z \int_0^{\xi_1} \cdots \int_0^{\xi_{n-2}} f^{(n)}_r(\xi_{n-1}) \, d\xi_{n-1} \cdots d\xi_1 + \sum_{j=0}^{n-2} \frac{f^{(j+1)}_r(0)}{j!} z^j \\
&= \int_0^z \int_0^{\xi_1} \cdots \int_0^{\xi_{n-2}} \Bigg( B_n(\xi_{n-1}) - \sum_{j=0}^{n-1} B_j(\xi_{n-1}) f^{(j)}_r(\xi_{n-1}) \Bigg) \, d\xi_{n-1} \cdots d\xi_1 \\
&\quad + \sum_{j=0}^{n-2} \frac{f^{(j+1)}_r(0)}{j!} z^j.
\end{align*}
Then by Lemma \ref{lem:qx} yields
\begin{align*}
&\sup_{a \in \mathbb{D}} \frac{(1 - |a|^2)^2}{K(1 - |a|^2)} \int_{\mathbb{D}} |f_r'(z)|^p (1 - |z|^2)^{p-4+s} (1 - |\sigma_a(z)|^2)^2 dA(z)
\\&= \sup_{a \in \mathbb{D}} \frac{(1 - |a|^2)^2}{K(1 - |a|^2)} \int_{\mathbb{D}} \left| \int_0^z \int_0^{\xi_1} \cdots \int_0^{\xi_{n-2}} \left( B_n(\xi_{n-1}) - \sum_{j=0}^{n-1} B_j(\xi_{n-1}) f_r^{(j)}(\xi_{n-1}) \right) d\xi_{n-1} \cdots d\xi_1 \right. \\
&\quad \left. + \sum_{j=0}^{n-2} \frac{f_r^{(j+1)}(0)}{j!} z^j \right|^p (1 - |z|^2)^{p-4+s} (1 - |\sigma_a(z)|^2)^2 dA(z) \\
&\leq \sup_{a \in \mathbb{D}} \frac{(1 - |a|^2)^2}{K(1 - |a|^2)} \int_{\mathbb{D}} \left| \int_0^z \int_0^{\xi_1} \cdots \int_0^{\xi_{n-2}} B_0(\xi_{n-1})f_r(\xi_{n-1}) d\xi_{n-1} \cdots d\xi_1 \right|^p(1 - |z|^2)^{p-4+s}
\\&\quad \cdot (1 - |\sigma_a(z)|^2)^2 dA(z) \\
&\quad + \sup_{a \in \mathbb{D}} \frac{(1 - |a|^2)^2}{K(1 - |a|^2)} \int_{\mathbb{D}} \left| \int_0^z \int_0^{\xi_1} \cdots \int_0^{\xi_{n-2}} \sum_{j=1}^{n-1} B_j(\xi_{n-1}) f_r^{(j)}(\xi_{n-1}) d\xi_{n-1} \cdots d\xi_1 \right|^p 
\\&\qquad \cdot (1 - |z|^2)^{p-4+s} (1 - |\sigma_a(z)|^2)^2 dA(z) 
\\&\quad+\sup_{a \in \mathbb{D}} \frac{(1 - |a|^2)^2}{K(1 - |a|^2)} \int_{\mathbb{D}} \left| \int_0^z \int_0^{\xi_1} \cdots \int_0^{\xi_{n-2}} B_n(\xi_{n-1}) d\xi_{n-1} \cdots d\xi_1 \right|^p (1 - |z|^2)^{p-4+s} 
\\&\quad\text{ }\text{ }\text{ }\cdot(1 - |\sigma_a(z)|^2)^2 dA(z)
\\&\quad + \sup_{a \in \mathbb{D}} \frac{(1 - |a|^2)^2}{K(1 - |a|^2)} \int_{\mathbb{D}} \left| \sum_{j=0}^{n-2} \frac{f_r^{(j+1)}(0)}{j!} z^j \right|^p (1 - |z|^2)^{p-4+s} (1 - |\sigma_a(z)|^2)^2 dA(z)
\\&:= N_{0,1} + N_{1,1} + N_2+N_3.
\end{align*}
Using lemma \ref{lem:jw}, we get
\begin{align*}
  N_{0,1}&\lesssim\sup_{a \in \mathbb{D}} \frac{(1 - |a|^2)^2}{K(1 - |a|^2)} \int_D \left( \int_0^z \int_0^{\xi_1} \cdots \int_0^{\xi_{n-2}} |B_0(\xi_{n-1})||f_r(\xi_{n-1})| \, d\xi_{n-1} \cdots d\xi_1 \right)^p 
\\&\quad\cdot (1 - |z|^2)^{p-4+s} (1 - |\sigma_a(z)|^2)^2 dA(z)
\\&\leq \|f_r\|_{B_{p}^{K}(s)}^{p} \sup_{a \in \mathbb{D}} \frac{(1 - |a|^2)^2}{K(1 - |a|^2)}  \int_D \bigg( \int_0^z \int_0^{\xi_1} \cdots \int_0^{\xi_{n-2}} \frac{|B_0(\xi_{n-1})|^{p}}{(1 - |\xi_{n-1}|^2)^{s}} \,  
 \\&\quad\cdot K(1-|\xi_{n-1}|^2)d\xi_{n-1} \cdots d\xi_1 \bigg)(1 - |z|^2)^{p-4+s} (1 - |\sigma_a(z)|^2)^2 dA(z)
\\&:= \|f_r\|_{B_{p}^{K}(s)}^{p} N_0.
\end{align*}
Using lemma \ref{lem:cl}, we estimate $N_{1,1}$,
\begin{align*}
N_{1,1} &= \sup_{a \in \mathbb{D}} \frac{(1 - |a|^2)^2}{K(1 - |a|^2)}\int_{\mathbb{D}} \bigg| \int_0^z \int_0^{\xi_1} \cdots \int_0^{\xi_{n-2}} \sum_{j=1}^{n-1} \sum_{i=0}^{j-1} (-1)^i \binom{j-1}{i} \bigg(B_j^{(i)}(\xi_{n-1}) 
\\&\quad\cdot f'_r(\xi_{n-1})\bigg)^{(j-i-1)} d\xi_{n-1} \cdots d\xi_1 \bigg|^p  (1 - |z|^2)^{p-4+s} (1 - |\sigma_a(z)|^2)^2 dA(z) 
\\&=  \sup_{a \in \mathbb{D}} \frac{(1 - |a|^2)^2}{K(1 - |a|^2)}\int_{\mathbb{D}} \bigg| \sum_{j=1}^{n-1} \sum_{i=0}^{j-1} (-1)^i \binom{j-1}{i} \int_{0}^{z} \int_{0}^{\xi_1} \cdots \int_{0}^{\xi_{n-2-(j-i-1)}} \bigg [\bigg(B_j^{(i)} 
\\&\quad \cdot f_r'\bigg)(\xi_{n-1-(j-i-1)}) - \sum_{t=0}^{j-i-2} \frac{\bigg( B_j^{(i)} f_r' \bigg)^{(t)} (0)}{t!} \xi_{n-1-(j-i-1)}^t \Bigg] d\xi_{n-1-(j-i-1)} \cdots d\xi_1 \Bigg|^p 
\\&\quad \cdot(1 - |z|^2)^{p-4+s} (1 - |\sigma_a(z)|^2)^2 dA(z) dA(z)
\\&\lesssim \sup_{a \in \mathbb{D}} \frac{(1 - |a|^2)^2}{K(1 - |a|^2)}\int_{\mathbb{D}} \bigg| \sum_{m=1}^{n-1} \int_0^z \int_0^{\xi_1} \cdots \int_0^{\xi_{m-1}} \sum_{k=1}^{m} (-1)^{m-k} \binom{n-k-1}{m-k} B_{n-k}^{(m-k)}(\xi_m) 
\\&\quad \cdot f'_r(\xi_m) d\xi_m \cdots d\xi_1 \bigg|^p(1 - |z|^2)^{p-4+s} (1 - |\sigma_a(z)|^2)^2 dA(z) 
\\&\quad+\sup_{a \in \mathbb{D}} \frac{(1 - |a|^2)^2}{K(1 - |a|^2)} \int_{\mathbb{D}} \left| \sum_{j=1}^{k-1} \sum_{i=0}^{j-1} \sum_{t=0}^{j-i-2} (-1)^i \binom{j-1}{i} \frac{(B_j^{(i)} f'_r)^{(t)}(0)}{(t+n+i-j)!} z^{t+n+i-j} \right|^p 
\\&\quad\quad\cdot (1 - |z|^2)^{p-4+s} (1 - |\sigma_a(z)|^2)^2 dA(z) 
\\&:= N_{111} + N_{112}.
\end{align*}
Using lemma \ref{lem:jw}, we have
\begin{align*}
N_{111} &\lesssim \sup_{a \in \mathbb{D}} \frac{(1 - |a|^2)^2}{K(1 - |a|^2)} \int_{\mathbb{D}} \left( \sum_{m=1}^{n-1} \int_0^z \int_0^{\xi_1} \cdots \int_0^{\xi_{m-1}} \left| \sum_{k=1}^m B_{n-k}^{(m-k)}(\xi_m) f_r'(\xi_m) \right| d\xi_m \cdots d\xi_1 \right)^p 
\\&\quad\cdot (1 - |z|^2)^{p-4+s} (1 - |\sigma_a(z)|^2)^2 dA(z) dA(z)
\\&\lesssim \|f_{r}\|^{p}_{B^K_p(s)} \sup_{a \in \mathbb{D}} \frac{(1 - |a|^2)^2}{K(1 - |a|^2)}
\int_{\mathbb{D}} \Biggl(\sum_{m=1}^{n-1} \int_0^z \int_0^{\xi_1} \cdots \int_0^{\xi_{m-1}} 
\frac{\bigl|\sum_{k=1}^m B^{(m-k)}_{n-k}(\xi_m)\bigr|^p}{(1-|\xi_m|^2)^{p+s}} 
\, 
\\&\quad \cdot K(1 - |\xi_m|^2)d\xi_m \cdots d\xi_1 \Biggr)(1-|z|^2)^{p-4+s} (1-|\sigma_a(z)|^2)^2 \, dA(z) \\
&= \|f_{r}\|^{p}_{B^K_p(s)} N_1.
\end{align*}
It is easy to check $N_{112}<\infty$ and $N_{3}<\infty$. Hence, the conditions of Theorem \ref{thm:fly} yield
$$\|f_{r}\|_{B_{p}^{K}(s)}\lesssim\frac{N_{112}+N_{3}+N_{2}}{1-N_{0}-N_{1}}<\infty,\quad0<r<1.$$
This implies $f\in B_{p}^{p}(s)$ by letting $r\rightarrow1^{-}$.
\qed

%%~~~~~~~~~~~~~~~~~~~~~~~~~~~~~~~~~~~~~~~~~~~~~~~~~reference~~~~~~~~~~~~~~~~~~~~~~~~~~~~~~~~~~~~~~~~~~~~~~~~~~~~~~~~~~~~~~~~~~~~~~~
%%~~~~~~~~~~~~~~~~~~~~~~~~~~~~~~~~~~~~~~~~~~~~~~~~~~~~~~~~~~~~~~~~~~~~~~~~~~~~~~~~~~~~~~~~~~~~~~~~~~~~~~~~~~~~~~

\end{document}